\documentclass[12pt]{article}
\usepackage{amsmath, amssymb, eucal, latexsym}
\usepackage[cp1251]{inputenc}
\usepackage[english]{babel}
\usepackage{graphicx}

\def\R{\Bbb R}

\def\a{\alpha}
\def\la{\lambda}
\def\D{\Delta}

\def\E{\mathsf {E}}

\newtheorem{theorem}{Theorem}
\newtheorem{proposition}[theorem]{Proposition}
\newtheorem{lemma}[theorem]{Lemma}
\newtheorem{cor}[theorem]{Corollary}
\newtheorem{note}[theorem]{Remark}

\begin{document}

 \begin{center}
 \textbf{ON SUM SETS OF SETS, HAVING SMALL PRODUCT SET}
 \end{center}

 \begin{center}
                                                         S. V. KONYAGIN, I. D. SHKREDOV
\footnote{
This work is supported by the RSF under a grant 14-50-00005.
}\\

    \end{center}

\bigskip

\begin{center}
    Annotation.
\end{center}

{\it \small
    We improve a result of Solymosi on sum--products in $\R$, namely, we prove that
    $\max{\{|A+A|,|AA|\}}\gg |A|^{\frac{4}{3}+c}$, where $c>0$ is an absolute constant.
    New lower bounds for sums of sets with small product set are found.
    Previous results are improved effectively
    for sets $A\subset \R$ with $|AA| \le |A|^{4/3}$.
}

\bigskip
\section{Introduction}
\bigskip

Let  $A,B\subset \R$ be finite sets.
Define the  \textit{sum set}, the \textit{product set} and \textit{quotient set} of $A$ and $B$ as
$$A+B:=\{a+b ~:~ a\in{A},\,b\in{B}\}\,,$$
$$AB:=\{ab ~:~ a\in{A},\,b\in{B}\}\,,$$
and
$$A/B:=\{a/b ~:~ a\in{A},\,b\in{B},\,b\neq0\}\,,$$
correspondingly.
The Erd\"{o}s--Szemer\'{e}di  conjecture \cite{ES} says that for any  $\epsilon>0$ one has
$$\max{\{|A+A|,|AA|\}}\gg{|A|^{2-\epsilon}} \,.$$
Roughly speaking, it asserts that an arbitrary subset of real numbers (or integers)
cannot has good additive and multiplicative structure, simultaneously.
At the moment the best result in this  direction is due to Solymosi \cite{soly}.

\begin{theorem}
    Let  $A\subset \R$ be a set.
    Then
\begin{equation}\label{f:Solymosi}
    |A+A|^2 |A/A|\,, \quad |A+A|^2 |AA|   \ge \frac{|A|^4}{4 \lceil \log |A| \rceil} \,.
\end{equation}
    In particular
\begin{equation}\label{f:Solymosi_max}
    \max{\{|A+A|,|AA|\}}\gg \frac{|A|^{4/3}}{\log^{1/3} |A|} \,.
\end{equation}
\label{t:Solymosi}
\end{theorem}
Here and below we suppose that $|A|\ge2$.

     It is easy to see that bound (\ref{f:Solymosi}) is tight up to logarithmic factors if
     the size of $A+A$ is small relatively  to $A$.
     The first part of the paper concerns the case where  the product $AA$ is small.
    We will write $a \lesssim b$ or $b \gtrsim a$ if $a = O(b \cdot \log^c |A|)$, $c>0$.
    In these terms inequality (\ref{f:Solymosi}) implies the following.

\begin{cor}
    Let $A\subset \R$ be a finite set and $K\ge 1$ be a real number.
    Suppose that $|A/A| \le K|A|$ or $|AA| \le K|A|$.
    Then
\begin{equation}\label{f:sol}
    |A+A| \gtrsim |A|^{\frac{3}{2}} K^{-\frac{1}{2}} \,.
\end{equation}
\label{c:sol}
\end{cor}

Estimate (\ref{f:sol}) was improved for small $K$, see e.g. references in paper \cite{s_sumsets}
(sharper bounds for {\it difference} of two sets, having small multiplicative doubling can be found in \cite{Sh_ineq}).
Here we give a result from
\cite{s_sumsets}.

\begin{theorem}
    Let $A\subset \R$ be a finite set and $K\ge 1$ be a real number.
    Suppose that $|A/A| \le K|A|$ or $|AA| \le K|A|$.
    Then
\begin{equation}\label{f:previous}
    |A+A| \gtrsim |A|^{\frac{58}{37}} K^{-\frac{42}{37}} \,.
\end{equation}
\label{t:previous}
\end{theorem}

It is easy to check that the bound of Theorem \ref{t:previous} is better than Corollary \ref{c:sol} for
$K\lesssim |A|^{\frac{5}{47}}$.

Let us formulate the first result of the article
(its  refined version is contained in  Theorem \ref{t:small_md} and Theorem \ref{t:small2_md} below).

\begin{theorem}
    Let $A\subset \R$ be a finite set and $K\ge 1$ be a real number.
Suppose that $|A/A| \le K|A|$ or $|AA| \le K|A|$.
    Then
\begin{equation*}\label{f:main}
    |A+A| \gtrsim |A|^{\frac{19}{12}} K^{-\frac{5}{6}}
\end{equation*}
    and
\begin{equation*}\label{}
 |A+A| \gtrsim |A|^{\frac{49}{32}} K^{-\frac{19}{32}} \,.
\end{equation*}
\label{t:main}
\end{theorem}

Theorem \ref{t:main} is stronger than  Theorem \ref{t:previous} and
refines estimate (\ref{f:sol}) for $K \lesssim |A|^{1/3}$.

In  Theorem  \ref{t:Sol+} we improve
bound (\ref{f:Solymosi_max}).

\begin{theorem}
    Let  $A\subset \R$ be a set.
    Then
$$\max{\{|A+A|,|AA|\}}\gg |A|^{\frac{4}{3}+c} \,,$$
where $c>0$ is an absolute constant.
\label{t:Sol_new}
\end{theorem}

Besides, a "critical"\, case of Solymosi's theorem, i.e.
the situation where the reverse inequality to (\ref{f:Solymosi}) takes place
is considered in the paper, see Proposition \ref{p:E(AA)_E(A)}.

We use a combination of methods from \cite{soly} and \cite{SS1} in our arguments.

We thank M.Z. Garaev for his advice to use a variant of
Balog--Szemer\'{e}di--Gowers theorem from \cite{BG}
and indication on inaccuracies in the previous version of the paper.

\bigskip
\section{Definitions and preliminary results}
\bigskip

The {\it additive energy $\E^{+} (A,B)$} between two sets $A$ and $B$ is the number of the solutions of the equation
(see \cite{TV})
$$
    \E^{+} (A,B) = |\{ a_1+b_1 = a_2+b_2 ~:~ a_1,a_2\in A\,, b_1,b_2\in B \}| \,.
$$
The {\it multiplicative energy $\E^{\times} (A,B)$} between two sets $A$ and $B$ is  the number of the solutions
of the equation (see \cite{TV})
$$
    \E^{\times} (A,B) = |\{ a_1 b_1 = a_2 b_2 ~:~ a_1,a_2\in A\,, b_1,b_2\in B \}| \,.
$$
In the case $A=B$ we write $\E^{+} (A)$ for $\E^{+} (A,A)$ and $\E^{\times} (A)$ for $\E^{\times} (A,A)$.
Having $\la\in A/A$, we  put $A_\la = A\cap \la A$. Clearly, if $0\not\in A$ then
\begin{equation}\label{f:energy_basic}
    \E^{\times} (A) = \sum_{\la \in A/A} |A_\la|^2
\end{equation}
and, similarly, for the energy $\E^{+} (A)$.
Finially, the Cauchy--Schwarz inequality implies for $0\not\in A$, $A_1\subset A$, $A_2\subset A$ that
\begin{equation}\label{f:energy_KB_subs}
    \E^{\times} (A_1,A_2) |A/A| \ge |A_1|^2|A_2|^2 \,, \quad \quad
\E^{\times} (A_1,A_2) |AA| \ge |A_1|^2|A_2|^2 \,.
\end{equation}
In particular
\begin{equation}\label{f:energy_KB}
    \E^{\times} (A) |A/A| \ge |A|^4 \,, \quad \quad \E^{\times} (A) |AA| \ge |A|^4 \,.
\end{equation}

Solymosi's Theorem \ref{t:Solymosi} can be derived from a slightly delicate result on an upper bound for the multiplicative energy of a set via its sum set, see \cite{soly}.
Estimation of the cardinality of the set
from the left hand side of (\ref{f:Solymosi-}) is the main task of our crucial
Lemma \ref{l:main_lemma}.

\begin{theorem}
    Let  $A,B\subseteq \R$ be a finite sets with $\min\{|A|, |B| \}\ge2$ and $\tau \ge 1$ be a real number. Then
\begin{equation}\label{f:Solymosi-}
    |\{ x ~:~ |A \cap xB| \ge \tau \}| \ll \frac{|A+A| |B+B|}{\tau^2} \,.
\end{equation}
    In particular
\begin{equation}\label{f:Solymosi_E}
    \E^{\times} (A,B) \ll |A+A| |B+B| \cdot \log (\min\{|A|, |B| \}) \,.
\end{equation}
\end{theorem}

We need in Lemma 7 from \cite{RR-NS}.
In paper \cite{SS1}, see Lemma 27, the same result was obtained with the additional factor $\log^2 d(A)$.

\begin{lemma}
    Let  $A\subset \R$ be a finite set.
    Then for any finite set $B\subset \R$ and an arbitrary real number $\tau \ge 1$ one has
\begin{equation}\label{f:d(A)}
    |\{ x\in A+B ~:~ |A\cap (x-B)| \ge \tau \}|
        \ll
        d(A)
    \cdot \frac{|A| |B|^2}{\tau^{3}} \,,
\end{equation}
    where
$$
    d(A) := \min_{C\neq \emptyset} \frac{|AC|^2}{|A||C|} \,.
$$
\label{l:d(A)}
\end{lemma}

Obviously, if
$|A/A| \le K|A|$ or $|AA| \le K|A|$ then  $d(A)$ does not exceed $K^2$.
The quantity $d(A)$ is a more delicate characteristic of a set than $|A/A|/|A|$ or $|AA|/|A|$.
For example, rough estimate (\ref{f:previous}) can be derived from a stronger one
\begin{equation}\label{f:previous_real}
    |A+A| \gtrsim |A|^{\frac{58}{37}} d(A)^{-\frac{21}{37}}
    \,,
\end{equation}
see \cite{s_sumsets}.

\bigskip

Lemma \ref{l:d(A)} implies the following  result.

\begin{cor}
    Let  $A_1,A_2,A_3 \subset \R$ be any finite sets
    and $\a_1,\a_2,\a_3$ be arbitrary nonzero numbers.
    Then the number of the solutions of the equation
\begin{equation}\label{f:gen_sigma}
    \sigma (\a_1 A_1, \a_2 A_2, \a_3 A_3) :=
    |\{ \a_1 a_1 + \a_2 a_2 + \a_3 a_3 = 0 ~:~ a_1\in A_1\,, a_2\in A_2\,, a_3\in A_3 \}|
\end{equation}
    does not exceed
    $O(d^{1/3} (A_1) |A_1|^{1/3} |A_2|^{2/3} |A_3|^{2/3})$.
\label{c:gen_sigma}
\end{cor}
{\bf Proof of the corollary.}
Without loosing of generality, we can suppose that $\a_1=1$.
Then the number of the solutions of equation (\ref{f:gen_sigma}) is
\begin{equation}\label{tmp:10.02.2015_1}
    \sigma:= \sum_{x\in (-\a_3 A_3)} |A_1\cap (x-\a_2 A_2)| \,.
\end{equation}
Let us arrange the values of $|A_1\cap (x-\a_2 A_2)|$ in decreasing order, that is
$|A_1\cap (x_1-\a_2 A_2)| \ge |A_1\cap (x_2 -\a_2 A_2)| \ge \dots$.
Using Lemma \ref{l:d(A)}, we obtain $|A_1\cap (x_j-\a_2 A_2)| \ll d^{1/3} (A_1) |A_1|^{1/3} |A_2|^{2/3} j^{-1/3}$.
Substitutioning the last bound in (\ref{tmp:10.02.2015_1}), we get
$$
    \sigma \ll d^{1/3} (A_1) |A_1|^{1/3} |A_2|^{2/3} |A_3|^{2/3}
$$
as required.

\bigskip

The last result of the section connects the quantity  $\E^{+} (A)$ with $|A/A|$ and $|AA|$.
We follow the arguments from \cite{ER} in the proof.

\begin{theorem}
    Let $A\subset \R$ be a finite set.
    Then
\begin{equation}\label{f:Sol_new_real}
    |A/A| |A|^{10}\log |A| \gg (\E^{+} (A))^4,\quad |AA| |A|^{10}\log |A| \gg (\E^{+} (A))^4 \,.
\end{equation}
\label{t:ER_new}
\end{theorem}
{\bf Proof of the theorem.} Without loosing of generality, we can suppose that
all elements of $A$ are positive.
For $x\in\R$ put
$$N(x)=|A\cap(x-A)| \,.$$
We have
\begin{equation}\label{1st2nd_mom}
\sum_{x\in A+A} N(x)=|A|^2,\quad \sum_{x\in A+A} N^2(x)= \E^+ (A) \,.
\end{equation}
Let
$$F=\left\{x\in A+A:\, N(x)>\frac{\E^+ (A)}{2|A|^2}\right\} \,.$$
Then
$$\sum_{x\not\in F}N^2(x)\le \sum_{x\not\in F}N(x) \cdot \frac{\E^+(A)}{2|A|^2} \,.$$
Using this and the first formula of (\ref{1st2nd_mom}), we obtain
$$\sum_{x\not\in F}N^2(x)\le |A|^2\frac{\E^+(A)}{2|A|^2}=\frac{\E^+(A)}{2} \,.$$
Applying (\ref{1st2nd_mom}) once more time, we get
\begin{equation}\label{sum_F}
\sum_{x\in F}N^2(x)\ge \frac{\E^+(A)}{2} \,.
\end{equation}
Put
$$U=\sum_{x\in F}N(x) \,.$$
Because of (\ref{sum_F}) and a trivial bound $N(x)\le|A|$, we have
\begin{equation}\label{est_U}
U \ge \frac{\E^+(A)}{2|A|} \,.
\end{equation}
Further, by the definition of the set  $F$
$$|F|\le\frac{2|A|^2 U}{\E^+(A)} \,.$$
Using this and inequality (\ref{est_U}), we obtain
\begin{equation}\label{est_F+A}
|F|+|A|\le\frac{4|A|^2 U}{\E^+(A)} \,.
\end{equation}
Let us consider the set of points from $\R^2$:
$$P=(A\cup F)\times(A\cup F)$$
and let us estimate  the number of collinear triples $T$  from $P$ (points in a triple are not necessarily distinct).
On the one hand, a general upper bound for the number of such triples in Cartesian products (\cite{TV}, Corollary~8.9) gives us
$$T\ll |A\cup F|^4\log|A| \,.$$
Because of (\ref{est_F+A}), it implies
\begin{equation}\label{tripple_upp}
T\ll\frac{|A|^8 U^4\log|A|}{(\E^+(A))^4} \,.
\end{equation}
On the other hand, for $x\in A$ put
$$F(x)=\{y\in A:\,x+y\in F\} \,.$$
Fixing $e, f\in A$, we have by (\ref{f:energy_KB_subs}) that there are at least
$$T(e,f)=F^2(e)F^2(f)/\min\{ |AA|,|A/A| \}$$
quadruples  $(a,b,c,d)$ such that $ab=cd$, $a,c\in F(e)$, $b,d\in F(f)$.
They form at least $T(e,f)$ collinear triples
$$(e,f),\,(e+a,f+d), (e+c,f+b) \,.$$
It follows that
$$T\ge\min(|AA|,|A/A|)^{-1}\sum_{e,f\in A}F^2(e)F^2(f)
=\min(|AA|,|A/A|)^{-1}\left(\sum_{e\in A}F^2(e)\right)^2 \,.$$
By the Cauchy--Schwarz inequality
$$\sum_{e\in A}F^2(e)\ge\left(\sum_{e\in A}F(e)\right)^2|A|^{-1}=U^2|A|^{-1} \,.$$
Whence
\begin{equation}\label{tripple_low}
T\ge\min(|AA|,|A/A|)^{-1}U^4|A|^{-2} \,.
\end{equation}
Combining estimates (\ref{tripple_upp}) and (\ref{tripple_low}), we obtain the required result.

\bigskip


\bigskip
\section{The proof of the main results}
\bigskip

We begin with a technical lemma.

Let $A\subset \R$, $0\not\in A$ be a finite set and $\tau>0$ be a real number.
Let also $S'_\tau$ be a set
$$
 S_\tau'\subset S_\tau := \{ \la ~:~ \tau < |A_\la| \le 2 \tau \} \subseteq A/A
$$
and for any nonzero $\alpha_1,\alpha_2,\alpha_3$ and different $\la_1,\la_2,\la_3\in S_\tau'$ one has
$$
        \sigma (\a_1 A_{\la_1}, \a_2 A_{\la_2}, \a_3 A_{\la_3})\le\sigma \,.
$$

\begin{lemma}
  Let $A\subset \R$, $0\not\in A$ be a finite set, $\tau>0$ be a real number,
    \begin{equation}\label{cond:main_lemma}
        32\sigma \le \tau^2 \le |A+A| \sqrt{\sigma} \,,
    \end{equation}
    and  $S_\tau'$,
    $\sigma$ are defined above.
    Then
\begin{equation}\label{f:main_lemma}
    |A+A|^2 \ge
     \frac{\tau^3 |S_\tau'|}{128 \sqrt{\sigma}}
     \,.
\end{equation}
\label{l:main_lemma}
\end{lemma}
{\bf Proof of the lemma.}
We follow the arguments from \cite{soly}.
Without loosing of generality, one can suppose that $A \subset \R^{+}$.
Consider the Cartesian product
$A\times A$ and the lines $l_\la$ of the form  $y=\la x$, where $\la \in A/A$.
Clearly, any line $l_\la$ intersects $A\times A$ under the points  $(x,\la x)$, $x\in A_\la$.
Put $\mathcal{A}_\la = l_\la \cap (A\times A)$.

Let $2 \le M \le |S_\tau'|$ be an integer parameter, which we will choose later.
Arrange the elements of the set $S_\tau'$ in increasing order 
and split it onto the groups of consecutive elements, each group has the size $M$.
We get  $k\ge \lfloor \frac{|S_\tau'|}{M} \rfloor \ge \frac{|S_\tau'|}{2M}$ such groups $U_j$.
Take the sets $\mathcal{A}_\la$
from
each of the group and consider all its sums.
Clearly, the sums belong  $(A+A)\times (A+A)$ and thus its total number does not exceed $|A+A|^2$.
On the other hand, by the inclusion--exclusion principle the number of such sums in any fixed group $U_j$ is at least
$$
    \rho_j := \tau^2 \binom{M}{2} - \sum_{\la_1,\dots,\la_4 \in U_j,\, \la_1 \neq \la_2,\,
\la_3 \neq \la_4,\,\{\lambda_1,\lambda_2\}\neq\{\lambda_3,\lambda_4\}}
        |\{ z ~:~ z\in (\mathcal{A}_{\la_1} + \mathcal{A}_{\la_2}) \cap (\mathcal{A}_{\la_3} + \mathcal{A}_{\la_4}) \}|
$$
\begin{equation}\label{f:rho_j}
            =
                \tau^2 \binom{M}{2} - \sum_{\la_1,\dots,\la_4 \in U_j,\, \la_1 \neq \la_2,\, \la_3 \neq \la_4,\,
\{\lambda_1,\lambda_2\}\neq\{\lambda_3,\lambda_4\}}
                    \mathcal{E} (\la_1,\dots,\la_4) \,.
\end{equation}
Fix  $\la_1,\dots,\la_4$ and prove that the quantity $\mathcal{E} (\la_1,\dots,\la_4)$
does not exceed $\sigma$.

Either all the numbers  $\la_1,\dots,\la_4$ are distinct or two of them coincide  but the other two are different
and differ from the first two numbers.
In any case there is a number, which differs from all of them.
Without loosing of generality, we can suppose that it is $\lambda_4$.
If
$$z=(z_1,z_2)\in (\mathcal{A}_{\la_1} + \mathcal{A}_{\la_2}) \cap (\mathcal{A}_{\la_3} + \mathcal{A}_{\la_4})$$
then $z_1=a_1+a_2=a_3+a_4$, $z_2=\la_1a_1+\la_2a_2=\la_3a_3+\la_4a_4$ for some
$a_j\in A_{\la_j}\,(j=1,2,3,4)$. It follows that
$$0=\la_1a_1+\la_2a_2-\la_3a_3-\la_4a_4-\la_4(a_1+a_2-a_3-a_4) $$
whence
$$(\la_1-\la_4)a_1+(\la_2-\la_4)a_2-(\la_3-\la_4)a_3=0 \,.$$
The number of tuples
$(a_1,a_2,a_3)$ satisfying the equation is
$$\sigma((\la_1-\la_4)A_{\la_1},(\la_2-\la_4)A_{\la_2},(\la_4-\la_3)A_{\la_3})\le\sigma \,.$$

Returning to formula (\ref{f:rho_j}) and using bound $\mathcal{E} (\la_1,\dots,\la_4)\le\sigma$,
we get
$$
    \rho_j \ge \tau^2 \binom{M}{2} - \sigma M^4 \,.
$$
Hence
$$
    |A+A|^2 \ge \frac{|S_\tau'|}{2M} \left(\tau^2 \binom{M}{2} - \sigma M^4  \right)
        \ge
            \frac{|S_\tau'|}{2M} \left(\frac{\tau^2 M^2}{4} - \sigma M^4 \right) \,.
$$
Put
$M = [\sqrt{\tau^2 / 8\sigma}]$.
The required inequality $M\ge2$ follows from the first condition
of (\ref{cond:main_lemma}).
Besides, if  we have $M \le |S_\tau'|$ then
$$
    |A+A|^2 \ge \frac{M \tau^2 |S_\tau'|}{16} \ge \frac{\tau^3 |S_\tau'|}{128 \sqrt{\sigma}}
$$
as required.
In contrary, suppose that $M > |S_\tau'|$ and assume that inequality  (\ref{f:main_lemma})
fails.
Then
$$|A+A|^2 < \frac{\tau^3 |S_\tau'|}{128 \sqrt{\sigma}} < \frac{\tau^3 M}{128 \sqrt{\sigma}}
< \frac{\tau^4}{256 \sigma} $$
with a contradiction to the RHS condition (\ref{cond:main_lemma}).
This concludes the proof of the lemma.

\bigskip

Let us prove the first part of Theorem \ref{t:main} which is our main result on sets with small product set.
 It is easy to see, that theorem below refines Solymosi's estimate (\ref{f:sol}) for $K \lesssim |A|^{1/4}$.

\begin{theorem}
    Let $A\subset \R$ be a finite set and $K\ge 1$ be a real number.
    Suppose that $|AA| \le K|A|$ or $|A/A| \le K|A|$.
    Then
\begin{equation}\label{f:small_md'}
    \E^\times (A) \ll K^{\frac{1}{4}}  |A|^{\frac{5}{8}} |A+A|^{\frac{3}{2}} (\log |A|)^{3/4} \,.
\end{equation}
    In particular
\begin{equation}\label{f:small_md}
    |A+A| \gg |A|^{\frac{19}{12}} K^{-\frac{5}{6}} (\log |A|)^{-1/2}\,.
\end{equation}
\label{t:small_md}
\end{theorem}
{\bf Proof of the theorem.}
 Estimate (\ref{f:small_md}) follows from (\ref{f:small_md'}) via inequality (\ref{f:energy_KB}) thus it is sufficient to prove (\ref{f:small_md'}).

Without loosing of generality, we can suppose that $0\notin A$.
Let $L=\log |A|$.
In the light of inequality  (\ref{f:Solymosi_E}) it is sufficient to check bound  (\ref{f:small_md'})
just for
$K^2 \le L^2 |A+A|^4 |A|^{-5}$.
From this bound and Solymosi's estimate (\ref{f:Solymosi}), we derive
\begin{equation}\label{tmp:12.02.2015_1}
    |A+A| \gg |A|^{\frac{11}{8}} L^{-\frac{1}{2}} \,.
\end{equation}
Further, because of $d(A)\le K^2$, we have
\begin{equation}\label{d(A)est}
d(A) \ll L^2 |A+A|^4 |A|^{-5}\,.
\end{equation}
Take a parameter
$\D = C L^{3/4} d^{1/8} (A) |A+A|^{3/2} |A|^{-11/8}$,
where $C>0$ is an absolute constant which we will choose later.
The constant $C$ depends on another constant $C_1>0$ which we will choose later as well.
By (\ref{d(A)est})
$$d(A)|A| \ll L^{3/2} d^{1/4}(A) |A+A|^3 |A|^{-11/4}$$
and we have for sufficiently large $C$ that
\begin{equation}\label{d(A)Deltaest}
C_1d(A) |A|\le \D^2\,.
\end{equation}

Further
\begin{equation}\label{tmp:06.03.2015_1}
    \E^{\times} (A) = \sum_{x} |A \cap xA|^2 \le \D |A|^2
        + \sum_{j\ge 1}\, \sum_{x ~:~ \D 2^{j-1} < |A \cap xA| \le \D 2^{j}} |A \cap xA|^2 \,.
\end{equation}
Let us note that
in formula (\ref{tmp:06.03.2015_1})
for
large enough $|A|$
it is sufficient to consider $j$ satisfying inequality
\begin{equation}\label{2to_j_est}
2^j \le |A|^{11/8} |A+A|^{-3/4}\,.
\end{equation}
Indeed, suppose in contrary that
$2^j > |A|^{11/8} |A+A|^{-3/4}$.
Then by inequality (\ref{tmp:12.02.2015_1}), we get
$$
    |A| \ge \D 2^{j} > C L^{3/4} d^{1/8} (A) |A+A|^{3/2}|A|^{11/8} |A|^{-11/8}  |A+A|^{-3/4}
        =
$$
$$
        =
            C L^{3/4} d^{1/8} (A) |A+A|^{3/4}
    \ge C L^{3/4} |A+A|^{3/4}  \gg_{C} |A|^{33/32}  L^{\frac{3}{8}}
$$
with a contradiction for large $|A|$.
Let $\tau = \D 2^{j-1}$ and $\sigma = \sigma (S_{\tau})$.
Take an arbitrary $\la \in S_\tau$.
By the definition of the set
$S_\tau$, we get  $d(A_\la) \le |A| \tau^{-1} d(A)$.
Applying  Corollary \ref{c:gen_sigma} and using the definition of the set $S_\tau$ once more time, we get
for any nonzero numbers
$\a_1,\a_2,\a_3$
    $$\sigma (\a_1 A_\la, \a_2 A_\la, \a_3 A_\la)\le\sigma,$$
where
\begin{equation}\label{tmp:16.02.2015_1}
    \sigma
        \ll
    (|A| \tau^{-1} d(A))^{1/3} \tau^{5/3}
\end{equation}
and we can take
$\sigma = M d^{1/3} (A) |A|^{1/3} \tau^{4/3}$, where $M>0$ is some constant.
Put $C_1=(32M)^{3}$, and the constant $C$ has chosen  such that inequality (\ref{d(A)Deltaest}) takes place.
It follows that
$$ \D^{2/3} \ge 32M d^{1/3} (A) |A|^{1/3}\,.$$
Hence for  $\tau\ge\Delta$, we have
$$ \tau^2 \ge 32M d^{1/3} (A) |A|^{1/3}\tau^{4/3} \,.$$
Thus the first condition of (\ref{cond:main_lemma}) takes place.

For any $j$ and sufficiently large $|A|$ in view of inequality  (\ref{2to_j_est}), we obtain
$$\tau =\D 2^{j-1} \le CL^{3/4} d^{1/8}(A) |A+A|^{3/4}
\le M^{3/8} |A|^{1/8} d^{1/8}(A) |A+A|^{3/4}\,.$$
It follows that
$$\tau^2 \le M^{1/2} |A|^{1/6} d^{1/6}(A) |A+A| \tau^{2/3} $$
and thus the second
inequality
of
(\ref{cond:main_lemma}) holds.

So, both conditions (\ref{cond:main_lemma}) for $\tau = \D 2^{j-1}$ take place
and we can apply inequality (\ref{f:main_lemma}) of the lemma to estimate
the cardinality of the set $S_{\D 2^{j-1}}$.
Using
(\ref{f:main_lemma}), (\ref{tmp:16.02.2015_1}), we get
$$
    \E^{\times} (A)
        \ll
            \D |A|^2 + \sum_{j\ge 1} \, \frac{d^{1/6} (A) |A|^{1/6} |A+A|^2}{2^{j/3} \D^{1/3}}
                \ll
                    \D |A|^2 \,.
$$
It follows that
$$
    \E^\times (A) \ll L^{3/4} d^{\frac{1}{8}} (A) |A|^{\frac{5}{8}} |A+A|^{\frac{3}{2}}
        \le
          L^{3/4} K^{\frac{1}{4}} |A|^{\frac{5}{8}} |A+A|^{\frac{3}{2}} \,.
$$
This completes the proof of the theorem.


\bigskip

In the next result we suppose that Solymosi's inequality (\ref{f:Solymosi})
cannot be improved.
We will show that the assumption implies lower bound for the additive energy of a set and
its product set $AA$.

\begin{lemma}
    Let $A\subset \R$, $0\not\in A$ be a finite set and $L\ge 1$ be a real number.
    Suppose that
\begin{equation}\label{cond:Sol_energyA/A}
    |A+A|^2 |A/A| \le L |A|^4 \,.
\end{equation}
    Then there is  $\tau \ge \E^{\times} (A)/(2|A|^2)$ and some sets $S'_\tau \subseteq S_\tau \subseteq A/A$,
    $|S_\tau| \tau^2 \gtrsim \E^{\times} (A)$, $|S'_\tau|\ge|S_\tau|/2$ such that
    for any element $\la$ from $S'_\tau$ one has
\begin{equation}\label{add_ener_low}
\E^{+} (A_\la) \gtrsim \tau^3 L^{-4}
\end{equation}
and
\begin{equation}\label{A/A_low}
|A_\la/A_\la| \gtrsim \tau^2 L^{-16}\,.
\end{equation}

Similarly, if
\begin{equation}\label{cond:Sol_energyAA}
    |A+A|^2 |AA| \le L |A|^4
\end{equation}
then there exists $\tau \ge \E^{\times} (A)/(2|A|^2)$ and some sets $S'_\tau \subseteq S_\tau \subseteq A/A$,
    $|S_\tau| \tau^2 \gtrsim \E^{\times} (A)$, $|S'_\tau|\ge|S_\tau|/2$ such that
    for any $\la \in S'_\tau$, we have (\ref{add_ener_low}) and
\begin{equation}\label{AA_low}
|A_\la A_\la| \gtrsim \tau^2 L^{-16}\,.
\end{equation}
\label{l:smallL}
\end{lemma}
{\bf Proof of the lemma.}
We consider the set
$A/A$ because the arguments in the case of the set $AA$ are similar.
One can assume
$$L=\max(1,|A+A|^2 |A/A| |A|^{-4})\,.$$
By Dirichlet principle there is
$\tau \ge \E^{\times} (A)/(2|A|^2)$
such that $|S_\tau| \tau^2 \gtrsim  \E^{\times} (A)$.
From (\ref{f:energy_KB}),
we have
\begin{equation}\label{tmp:energy_trivial}
    |S_\tau| \tau^2 \gtrsim \frac{|A|^4}{|A/A|}\,.
\end{equation}
If $|S_\tau| \ge 2$ then by  $S_\tau''$ denote the set of cardinality $[|S_\tau|/2]$ consisting all  $\la\in S_\tau$
with the minimal additive energy $\E^{+} (A_\la)$ and put $S_\tau'=S_\tau\setminus S_\tau''$.
It is sufficient to check that for some $\lambda\in S_\tau''$ one has
\begin{equation}\label{lower_add_ener}
\E^{+} (A_\la) \gtrsim \tau^3 L^{-4} \,.
\end{equation}
In the case $|S_\tau| = 1$ we put $S_\tau'=S_\tau''=S_\tau$
and it is sufficient to check inequality (\ref{lower_add_ener}) again.

Put
$\sigma := \max_{\la \in S_\tau''} \sqrt{2\tau \E^{+} (A_\la)}$.
Bound (\ref{add_ener_low}) follows from the inequality
\begin{equation}\label{lower_sigma}
\sigma \gtrsim \tau^2 L^{-2} \,,
\end{equation}
which is aim of our proof.

By the Cauchy--Schwarz inequality for any $\a,\beta \neq 0$ and arbitrary sets
$A_{\la_1}, A_{\la_2}, A_{\la_3}$, $\la_1,\la_2, \la_3 \in S_\tau''$ one has
$$
    \sigma(A_{\la_1}, \a A_{\la_2}, \beta A_{\la_3})
        \le
            |A_{\la_2}|^{1/2} (\E^{+} ( A_{\la_1}, \beta A_{\la_3}))^{1/2}
                \le
$$
$$
                \le
                    (2\tau)^{1/2} \E^{+} ( A_{\la_1})^{1/4} \E^{+} ( A_{\la_3})^{1/4}
                        \le
                            \sigma  \,.
$$

If both conditions (\ref{cond:main_lemma}) of Lemma \ref{l:main_lemma}
(with $S_\tau''$ instead of $S_\tau'$) take place then we have
$$
    |A+A|^2 \ge
     \frac{\tau^3 |S_\tau''|}{128 \sqrt{\sigma}}
 \ge \frac{\tau^3 |S_\tau|}{384 \sqrt{\sigma}}\,.
$$
Using condition (\ref{cond:Sol_energyA/A}), we get
\begin{equation}\label{tmp:12.02.2015_2}
    \sigma^{1/2} \gg \frac{|S_\tau| \tau^3 |A/A|}{|A|^4 L} \,.
\end{equation}
Substituting inequality (\ref{tmp:energy_trivial}) into  (\ref{tmp:12.02.2015_2}), we get  (\ref{lower_sigma}).

If the first condition (\ref{cond:main_lemma}) does not hold then we obtain (\ref{lower_sigma}) immediately.
Suppose that the second condition (\ref{cond:main_lemma})
fails,
that is
$\tau^2 > |A+A|\sqrt{\sigma}$.
By inequality (\ref{f:energy_KB}) for sums, we have a lower bound for $\sigma$, namely,
$\sigma^2 \ge 2 \tau^5 |A+A|^{-1}$.
But then
$$
    \tau^8 > |A+A|^4 \cdot 2\tau^5 |A+A|^{-1}
$$
with a contradiction, because, clearly, the parameter $\tau$ does not exceed the size of $A$.

Thus, we have proved inequality (\ref{add_ener_low}). Using Theorem  \ref{t:ER_new}, we obtain
inequality (\ref{A/A_low}). This concludes the proof of the lemma.

\bigskip

Now let us obtain the second main result of the paper, concerning the sets with small product set.
It is easy to see that we improve inequality (\ref{f:sol}) for $K \lesssim |A|^{1/3}$.

\begin{theorem}
    Let $A\subset \R$ be a finite set and $K\ge 1$ be a real number.
    Suppose that $|AA| \le K|A|$ or $|A/A| \le K|A|$.
    Then
\begin{equation}\label{f:small2_md}
    |A+A| \gtrsim |A|^{\frac{49}{32}} K^{-\frac{19}{32}} \,.
\end{equation}
\label{t:small2_md}
\end{theorem}
{\bf Proof of the theorem.}
Consider the situation  where
$|A/A| \le K|A|$. The case $|AA| \le K|A|$ is similar.
One can suppose that $0\notin A$.
Let us apply Lemma \ref{l:smallL}, where
$$L=\max(1,|A+A|^2 |A/A| |A|^{-4})\,.$$
Take any $\la$ from $S'_\tau$ and use inequality (\ref{A/A_low}) combining with the lower bound for $\tau$.
It gives us
$$|A/A| \ge |A_\la/A_\la| \gtrsim \tau^2 L^{-16} \ge (\E^{\times} (A))^2|A|^{-4} L^{-16} \,.$$
Further, because of  (\ref{f:energy_KB}), we have
$$|A/A| \gtrsim |A|^4 |A/A|^{-2}L^{-16}.$$
It follows that
$$L \gtrsim |A|^{1/4} |A/A|^{-3/16}.$$
After some simple calculations we obtain the result.

\bigskip

Theorem \ref{t:small2_md} improves Theorem \ref{t:small_md} for $K \gtrsim |A|^{5/23}$.

\bigskip

Let us obtain a result on multiplicative energies of $A/A$, $AA$ in "critical case".

\begin{proposition}
    Let $A\subset \R$ be a finite set.
    If  condition (\ref{cond:Sol_energyA/A}) takes place then
\begin{equation}\label{f:Sol_energy''}
    \E^{\times} (A/A) \gtrsim \frac{(\E^{\times} (A))^3}{L^{32} |A|^4}
    \,.
\end{equation}
    If condition (\ref{cond:Sol_energyAA}) holds  then
\begin{equation}\label{f:Sol_energy''_AA}
    \E^{\times} (AA) \gtrsim \frac{(\E^{\times} (A))^3}{L^{32} |A|^4}
    \,.
\end{equation}
\label{p:E(AA)_E(A)}
\end{proposition}
{\bf Proof of the proposition.}
Without loosing of generality, we can suppose that $0\notin A$.
Let us begin with inequality (\ref{f:Sol_energy''}).
Put $\Pi = A/A$.
Using Lemma \ref{l:smallL}, we find the number $\tau$ and the set $S'_\tau$ satisfying all implications of the lemma.
By the Katz--Koester inclusion (see \cite{kk}), namely  $A_\la / A_\la \subseteq \Pi \cap \la \Pi$,
we see that  for all $\la \in S'_\tau$
the following holds
$$
   |\Pi \cap \la \Pi| \ge   |A_\la / A_\la| \gtrsim \tau^2 L^{-16} \,.
$$
Hence
\begin{equation}\label{tmp:14.02.2015_1}
    \sum_{\la \in S'_\tau} |\Pi \cap \la \Pi| \gtrsim L^{-16} \tau^2 |S_\tau|
        \gtrsim
            L^{-16} \E^{\times} (A) \,.
\end{equation}
Using the last bound as well as the Cauchy--Schwarz inequality, we get  (\ref{f:Sol_energy''}).

Now put $\Pi' =AA$.
Then by the Katz--Koester inclusion, we have $A_\la  A_\la \subseteq \Pi' \cap \la \Pi'$ and the previous arguments can be applied.
This completes the proof of the proposition.

\bigskip

Thus, if
$|A/A| \lesssim |A|^{4/3}$ and $L \lesssim 1$ then inequality (\ref{f:Sol_energy''}) and bound (\ref{f:energy_KB}) imply
$\E^{\times} (A/A) \gtrsim L^{-32} |A/A|^3 \gtrsim |A/A|^3$.
In other words the multiplicative energy of the set $A/A$ is close to its maximal possible value.
We use the observation in the proof of the final
result
 of the paper.

\begin{theorem}
    Let  $A\subset \R$ be a set.
    Then for any
$c<\frac{1}{20598}$
    one has
\begin{equation}\label{f:Solymosi_max+}
    \max{\{|A+A|,|A/A|\}}\gg |A|^{\frac{4}{3}+c}
\end{equation}
    and
\begin{equation}\label{f:Solymosi_max+_AA}
    \max{\{|A+A|,|AA|\}}\gg |A|^{\frac{4}{3}+c}
    \,.
\end{equation}
\label{t:Sol+}
\end{theorem}
{\bf Proof of the theorem.}
We prove estimate (\ref{f:Solymosi_max+}) because inequality (\ref{f:Solymosi_max+_AA}) can be obtained  similarly.
Without loosing of generality, suppose that $0\notin A$.
Now assume  that inequality (\ref{cond:Sol_energyA/A}) holds with some parameter $L$.
Let also  $|A/A|^3 \le L' |A|^{4}$.
Our task is to find a lower bound for quantities $L$, $L'$.
Using Lemma \ref{l:smallL}, we have $\tau \ge \E^{\times} (A)/(2|A|^2)$ and sets
$S'_\tau \subseteq S_\tau \subseteq A/A$,
$|S_\tau| \tau^2 \gtrsim \E^{\times} (A)$, $|S'_\tau| \gtrsim |S_\tau|$ such that
for any element $\la$ from $S'_\tau$
one has
$|A_\la / A_\la| \gtrsim L^{-16} \tau^2$.
Using this as well as the Katz--Koester inclusion, we obtain
$$
    \sum_{x\in AA/AA} |S'_\tau \cap x (A/A)|
   = \sum_{\la \in S'_\tau} |A/A \cap \la (A/A)|
        \ge
            \sum_{\la \in S'_\tau} |A_\la / A_\la|
                \gtrsim L^{-16} \tau^2 |S_\tau| \,.
$$
In view of the last bound and the Cauchy--Schwarz inequality, we get
$$
    |A/A| \E^{\times} (S'_\tau, A/A) = |A/A| \sum_{x} |S'_\tau \cap x (A/A)|^2
        \gtrsim L^{-32} \tau^4 |S_\tau|^2 \,.
$$
Applying the Cauchy--Schwarz inequality once more time,
we obtain
$$
    \E^{\times} (S'_\tau) \gtrsim L^{-64} \tau^8 |S_\tau|^4 |A/A|^{-2} (\E^{\times} (A/A))^{-1}
        \gtrsim L^{-64} \E^{\times} (A) \tau^{6} |A/A|^{-5} |S_\tau|^3
            =
$$
$$
    =
    \eta |S_\tau|^3\,,
$$
where $\eta = L^{-64} \E^{\times} (A) \tau^{6} |A/A|^{-5}$.
We have
$$
\eta \gg L^{-64} \E^{\times} (A) \left(\E^{\times} (A) |A|^{-2}\right)^6 |A/A|^{-5}
= L^{-64} \E^{\times} (A)^7 |A|^{-12} |A/A|^{-5}
$$
$$
 \ge L^{-64} \left( |A|^4 |A/A|^{-1}\right)^7 |A|^{-12} |A/A|^{-5}
= L^{-64} |A|^{16} |A/A|^{-12} \ge L^{-64} (L')^{-4}\,.
$$
In other words
$$
    \E^{\times} (S'_\tau) \gtrsim L^{-64} (L')^{-4} |S_\tau|^3\,.
$$

By Balog--Szemer\'{e}di--Gowers Theorem \cite{BG} (see also \cite{S_BSzG}) there is
a set $S''_\tau \subseteq S'_\tau$,
$|S''_\tau| \gtrsim \eta |S_\tau|$ such that $|S''_\tau / S''_\tau| \lesssim \eta^{-4} |S''_\tau|^3 |S'_\tau|^{-2}$.
Because of  $S''_\tau \subseteq S_\tau$, we obtain
$$
    \sum_{a\in A} |A\cap a S''_\tau| = \sum_{\la \in S''_\tau} |A\cap \la A|
\gg \tau |S''_\tau|
$$
and hence there is $a\in A$ such that for the set  $A' := A\cap a S''_\tau$
one has
\begin{equation}\label{A'first_est}
|A'| \gg \tau |S''_\tau| |A|^{-1}\,.
\end{equation}
It follows that
$$
    d(A') \le \frac{|A' / S''_\tau|^2}{|A'| |S''_\tau|} \ll \frac{|S''_\tau / S''_\tau|^2 |A| }{\tau |S''_\tau|^2}
        \lesssim
            \eta^{-8} \frac{|A|}{\tau} \cdot \frac{|S''_\tau|^4}{|S_\tau|^4}\,.
$$
Using inequalities (\ref{f:energy_KB}), (\ref{f:previous_real}) and
the estimate for $d(A')$, we get
$$
    |A+A| \ge |A'+A'| \gtrsim |A'|^{\frac{58}{37}} d(A')^{-\frac{21}{37}}
        \gtrsim
            (\tau |S''_\tau| |A|^{-1})^{\frac{58}{37}} (\eta^{8} \tau |A|^{-1} |S_\tau|^4 |S''_\tau|^{-4})^{\frac{21}{37}}
$$
$$
    \gtrsim
        |S_\tau|^{\frac{58}{37}} (\tau |A|^{-1})^{\frac{79}{37}} \eta^{\frac{168}{37}}
            \gtrsim
                (\E^\times (A))^{\frac{58}{37}} |A|^{-\frac{79}{37}} \eta^{\frac{168}{37}} \tau^{-1}
$$
$$
    =
(\E^\times (A))^{\frac{58}{37}} |A|^{-\frac{79}{37}} \eta^{\frac{971}{222}}
             (L^{-64} \E^{\times} (A) |A/A|^{-5})^{\frac{1}{6}}
$$
$$
    \gtrsim
        (\E^\times (A))^{\frac{58}{37}} |A|^{-\frac{79}{37}} (L^{-64} (L')^{-4})^{\frac{971}{222}}
             (L^{-64} \E^{\times} (A) |A/A|^{-5})^{\frac{1}{6}}
  $$
$$
    =
        L^{-\frac{10752}{37}} (\E^{\times} (A))^{\frac{385}{222}} |A|^{-\frac{79}{37}}
            (L')^{-\frac{1942}{111}} |A/A|^{-\frac{5}{6}}
$$
$$
        \gtrsim
    L^{-\frac{10752}{37}} (|A|^4 |A/A|^{-1})^{\frac{385}{222}} |A|^{-\frac{79}{37}}
            (L')^{-\frac{1942}{111}} |A/A|^{-\frac{5}{6}}
$$
$$
    =
        L^{-\frac{10752}{37}} (L')^{-\frac{1942}{111}} |A|^{\frac{533}{111}} |A/A|^{-\frac{95}{37}}
        \gtrsim
        L^{-\frac{10752}{37}} (L')^{-\frac{1942}{111}} |A|^{\frac{533}{111}} ((L')^{1/3} |A|^{4/3})^{-\frac{95}{37}}
$$
$$
    =
        |A|^{\frac{51}{37}} L^{-\frac{10752}{37}} (L')^{-\frac{679}{37}} \,.
$$

The last estimate is greater than $|A|^{4/3}$ by some power of $|A|$.
Easy calculations show that one can take any number less than
$\frac{1}{20598}$ for the constant $c$.
This concludes the proof.


\begin{note}
    It seems likely
    that
    the arguments of the proof of Theorem \ref{t:Sol+}
    allow to improve slightly the lower bound for the size of  $A+A$ of
    Theorem \ref{t:small2_md} in the regime where $K \lesssim |A|^{1/3}$.
    We did not make such calculations.
\end{note}



{}
\bigskip
\newpage

\noindent{S.V.~Konyagin\\
Steklov Mathematical Institute,\\
ul. Gubkina, 8, Moscow, Russia, 119991}
and
\\
MSU,\\
Leninskie Gory, Moscow, Russia, 119992\\
{\tt konyagin@mi.ras.ru}

\bigskip

\noindent{I.D.~Shkredov\\
Steklov Mathematical Institute,\\
ul. Gubkina, 8, Moscow, Russia, 119991}
and
\\
IITP RAS,  \\
Bolshoy Karetny per. 19, Moscow, Russia, 127994\\
{\tt ilya.shkredov@gmail.com}


\begin{thebibliography}{99}


\bibitem{BG}
{\sc J.~Bourgain, M.Z.~Garaev, }
\emph{On a variant of sum-product estimates and explicit exponential sum bounds in prime fields, }
Math. Proc. Cambridge Philos. Soc. 146 (2009), no.1, 1--21.



\bibitem{ER}
{\sc G.~Elekes, I.~Ruzsa, }
\emph{Few sums, many products, }
Studia Sci. Math. Hungar. 40 (2003), no. 3, 301--308.


\bibitem{ES}
{\sc P.~Erd\H{o}s, E.~Szemer\'{e}di, }
\emph{On sums and products of integers, }
Studies in pure mathematics, 213--218, Birkh\"auser, Basel, 1983.


\bibitem{kk}
{\sc N.H.~Katz, P.~Koester, } {\em On additive doubling and energy,
} SIAM J. Discrete Math., 24 (2010), 1684--1693.




\bibitem{RR-NS}
{\sc O.E.~Raz, O.~Roche--Newton, M.~Sharir, }
{\em Sets with few distinct distances do not have heavy lines, }
arXiv:1410.1654v1 [math.CO] 7 Oct 2014.


\bibitem{S_BSzG}
{\sc T.~Schoen, }
{\em New bounds in Balog--Szemer\'{e}di--Gowers theorem, }
Combinatorica, {\bf 34}:5 (2014), 1--7.



\bibitem{SS1}
{\sc T. Schoen, I.D. Shkredov. }
{\em Higher moments of convolutions, }
J. Number Theory 133 (2013), no. 5, 1693--1737.



\bibitem{Sh_ineq}
{\sc I.D. Shkredov, }
{\em Some new inequalities in additive combinatorics, }
Moscow J. Combin. Number Theory 3 (2013), 237--288.


\bibitem{s_mixed}
{\sc I.D. Shkredov, }
{\em Some new results on higher energies, }
Transactions of MMS, 74:1 (2013), 35--73.


\bibitem{s_energy}
{\sc I.D. Shkredov, }
{\em Energies and structure of additive sets, }
Electronic Journal of Combinatorics, 21(3) (2014), \#P3.44, 1--53.


\bibitem{s_sumsets}
{\sc I.D. Shkredov, }
{\em On sums of Szemer\'{e}di--Trotter sets, }
arXiv:1410.5662v1 [math.CO] 21 Oct 2014.



\bibitem{soly}
{\sc J. Solymosi, }
{\em Bounding multiplicative energy by the sumset, }
Advances in Mathematics Volume 222, Issue 2, (2009), 402--408.




\bibitem{TV}
{\sc T. Tao, V. Vu, }
{\em Additive Combinatorics, }
Cambridge University Press (2006).




\end{thebibliography}
\end{document}